\documentclass{amsart}[12pt]
\usepackage{amsmath}
\usepackage{amstext}
\usepackage{amssymb}
\usepackage{amsthm}
\usepackage{amscd}

\usepackage{upref}

\swapnumbers
\theoremstyle{plain}
\newtheorem{thm}{Theorem}[section]
\newtheorem{lem}[thm]{Lemma}

\newtheorem{cor}[thm]{Corollary}

\theoremstyle{definition}

\theoremstyle{remark}

\DeclareMathOperator{\Coker}{Coker}

\DeclareMathOperator{\Ass}{Ass}

\def\nn{\relax\ifmmode{\mathbb N_{0}}\else$\mathbb N_{0}$\fi}
\def\lra{\longrightarrow}

\begin{document}

\title[An example of an infinite set of associated primes of a local cohomology module]
{An example of an infinite set of associated primes of a local cohomology module}
\author{MORDECHAI KATZMAN}
\address{Department of Pure Mathematics,
University of Sheffield, Hicks Building, Sheffield S3 7RH, United Kingdom\\
{\it Fax number}: 0044-114-222-3769}
\email{M.Katzman@sheffield.ac.uk}


\date{\today}



\maketitle

\setcounter{section}{-1}
\section{\bf Introduction}
\label{in}

Let $(R,m)$ be a local Noetherian ring, let $I\subset R$ be any ideal and let $M$ be a finitely generated $R$-module.
It has been long conjectured that the local cohomology modules $H^i_I(M)$ have finitely
many associated primes for all $i$ (see Conjecture 5.1 in \cite{H} and \cite{L}.)

If $R$ is not required to be local these sets of associated primes may be infinite, as shown by Anurag Singh
in \cite{S}, where he constructed an example of a
local cohomology module of a finitely generated module over  a finitely generated $\mathbb{Z}$-algebra with
infinitely many associated primes.
This local cohomology module has $p$-torsion for all primes $p\in \mathbb{Z}$.

However, the question of the finiteness of the set of associated primes of
local cohomology modules defined over local rings and over $k$-algebras (where $k$ is a field)
has remained open until now. In this paper I settle this question by constructing a
local cohomology module of a local finitely generated $k$-algebra
with an infinite set of associated primes, and I do this
for any field $k$.

\section{\bf The example}

Let $k$ be any field, let $R_0=k[x,y,s,t]$ and let
$S=R_0[u,v]$. Define a grading on $S$ by declaring $\deg(x)=\deg(y)=\deg(s)=\deg(t)=0$ and
$\deg(u)=\deg(v)=1$. Let $f=sx^2v^2 - (t+s) xy uv + t y^2 u^2$ and let $R=S/fS$.
Notice that $f$ is homogeneous and hence $R$ is graded. Let $S_+$ be the ideal of $S$
generated by $u$ and $v$ and
let $R_+$ be the ideal of $R$
generated by the images of $u$ and $v$.

Consider the local cohomology module $H^2_{R_+}(R)$:
it is homogeneously isomorphic to
$H^2_{S_+}(S/fS)$ and we can use the exact sequence
$$
H^2_{S_+}(S)(-2) \stackrel{f}{\lra} H^2_{S_+}(S) \lra
H^2_{S_+}(S/fS) \lra 0
$$
of graded $R$-modules and homogeneous homomorphisms (induced from
the exact sequence
$$
0 \lra S(-2) \stackrel{f}{\lra} S \lra S/fS \lra 0 )
$$
to study $H^2_{R_+}(R)$.
Furthermore, we can realize
$H^2_{S_+}(S)$ as the module $R_0[u^-,v^-]$ of inverse
polynomials described in \cite[12.4.1]{BS}: this graded $S$-module
vanishes beyond degree $-2$, and, for each $d \geq 2$, its $(-d)$-th component
is a free $R_0$-module of rank $d-1$ with base
$\left(u^{-\alpha}v^{-\beta}\right)_{\alpha,\beta>0,~\alpha + \beta = -d}.$
We will study the graded components of $H^2_{S_+}(S/fS)$ by considering
the cokernels of the $R_0$-homomorphisms
$$
f_{-d}: R_0[u^-,v^-]_{-d-2} \lra R_0[u^-,v^-]_{-d} \quad (d \geq 2)
$$
given by multiplication by $f$. In order to represent these
$R_0$-homomorphisms between free $R_0$-modules by matrices, we
specify an ordering for each of the above-mentioned bases by
declaring that
$$u^{\alpha_1}v^{\beta_1} < u^{\alpha_2} v^{\beta_2}$$
(where $\alpha_1, \beta_1, \alpha_2, \beta_2<0$ and $\alpha_1+\beta_1= \alpha_2+ \beta_2$) precisely when
$\alpha_1 > \alpha_2$.
If we use this ordering for both the source and target of each $f_d$, we can see that
each $f_d$ $(d\geq 2)$ is given by multiplication on the left by the tridiagonal
$d-1$ by $d+1$ matrix
$$
A_{d-1}:=
\left(
\begin{array}{ccccccc}
  s x^2 & -xy(t+s) & ty^2 & 0 &\dots  &  & 0 \\
  0 & s x^2 & -xy(t+s) & ty^2 & 0\dots &  & 0 \\
  0  & 0 & s x^2 & -xy(t+s) & ty^2 \dots&  & 0 \\
    &   & & \ddots & &  &  \\
  0 & \dots & & s x^2 & -xy(t+s) & & ty^2  \\
\end{array}
\right) .
$$
We also define
$$
\overline{A}_{d-1}:=
\left(
\begin{array}{ccccccc}
  s & -(t+s) & t & 0 &\dots  &  & 0 \\
  0 & s & -(t+s) & t & 0 \dots&  & 0 \\
  0 & 0 & s & -(t+s) & t\dots &  &   0 \\
    &   & & \ddots & &  &  \\
  0 & \dots & & s & -(t+s) & & t  \\
\end{array}
\right)
$$
obtained by substituting $x=y=1$ in $A_{d-1}$.

Let also  $\tau_i=(-1)^i (t^i+s t^{i-1}+\dots+s^{i-1} t+s^i)$.

\bigskip
\begin{lem}\hfil\label{lemma1}
\begin{itemize}
\item[(i)]
Let $B_i$ be the submatrix of
$\overline{A}_{i}$ obtained by deleting its first and last columns.
Then $\det B_i = \tau_i$ for all $i\geq 1$.
\item[(ii)] Let $\mathcal{S}$ be an infinite set of positive integers.
Suppose that either $k$ has characteristic zero or that $k$ has prime characteristic $p$ and
$\mathcal{S}$ contains infinitely many integers of the form $p^m-2$.
The ($k[s,t]$-)irreducible factors of $\left\{ \tau_i \right\}_{i\in \mathcal{S}}$
form an infinite set.
\end{itemize}
\end{lem}
\begin{proof}
We prove the first statement by induction on $i$. Since
$$\det B_1=\det \left(-t-s\right)=-t-s \textrm{\ and\ }
\det B_2=\det \left( \begin{array}{cc} -t-s & t \\ s & -t-s \end{array} \right)=t^2+st+s^2 ,$$
the lemma holds for $i=1$ and $i=2$. Assume now that $i\geq 3$.
Expanding the determinant of $B_i$ by its first row and applying the induction hypothesis we obtain
\begin{align*}
\det B_i &= (-t-s) \det B_{i-1} - s t \det B_{i-2}\\
&=(-1)^{i-1}(-t-s)(t^{i-1}+\dots+s^{i-2}t+s^{i-1})- (-1)^{i-2} st (t^{i-2}+\dots+s^{i-3}t+s^{i-2})\\
&=(-1)^{i} \left[ (t^{i}+\dots+s^{i-2} t^2+ s^{i-1}t) +
(s t^{i-1}+\dots+s^{i-1}t+s^i) - (s t^{i-1}+\dots+s^{i-2} t^2+ s^{i-1} t) \right]\\
&=(-1)^i (t^{i}+st^{i-1}+\dots+s^{i-1}t+s^i) .
\end{align*}

We now prove the second statement.
Define $\sigma_i=t^i+t^{i-1}+\dots+t+1$ and notice that it is enough to show that
the set of irreducible factors of
$\{ \sigma_i \}_{i\in \mathcal{S}}$ is infinite.
Let $\mathcal{I}$ be the set of irreducible factors of $\left\{ \sigma_i \right\}_{i\in \mathcal{S}}$.
If $k$ has characteristic zero consider
$\mathbb{Q}[\mathcal{I}] \supseteq \mathbb{Q}$, the splitting field of this set of irreducible factors.
If $\mathcal{I}$ is finite,
$\mathbb{Q}[\mathcal{I}] \supseteq \mathbb{Q}$ is finite extension
which contains all $i$th roots of unity for \emph{all} $i\in \mathcal{S}$, which is impossible.

Assume now that $k$ has prime characteristic $p$. Let $\mathbb{F}$ be the algebraic closure of
the prime field of $k$.
For any positive integer $m$
$$\frac{d}{dt} t(t^{p^m-1}-1) = -1$$
so $\sigma_{p^m-2}=(t^{p^m-1}-1)/(t-1)$ has $p^m-2$ distinct roots in $\mathbb{F}$ and, therefore,
the roots of $\{ \sigma_s \}_{s\in\mathcal{S}}$ form an infinite set.
\end{proof}

\bigskip
\begin{thm}\label{theorem1}
For every $d\geq 2$ the $R_0$-module $H^2_{R_+}(R)_{-d}$ has $\tau_{d-1}$-torsion.
Hence $H^2_{R_+}(R)$ has infinitely many associated primes.
\end{thm}
\begin{proof}
For the purpose of this proof we introduce a bigrading in
$R_0$ by declaring $\deg(x)=(1,0)$, $\deg(y)=(1,1)$ and $\deg(t)=\deg(s)=(0,0)$.

We also introduce a  bigrading on the free
$R_0$-modules $R_0^n$ by declaring
$\deg(x^\alpha y^\beta s^a t^b \mathbf{e}_j)=(\alpha+\beta, \beta+j)$
for all non-negative integers $\alpha, \beta, a,b $ and all $1 \leq j \leq n$. Notice that
$R_0^n$ is a bigraded $R_0$-module when $R_0$ is equipped with the bigrading mentioned above.

Consider the $R_0$-module $\Coker A_{d-1}$; the columns of $A_{d-1}$ are bihomogeneous of bidegrees
$$(2,1), (2,2), \dots, (2,d+1) .$$
We can now consider $\Coker A_{d-1}$ as a $k[s,t]$ module generated by the natural images of
$x^\alpha y^\beta \mathbf{e}_j$
for all non-negative integers $\alpha, \beta$ and all $1 \leq j \leq d-1$. The $k[s,t]$-module
of relations among these generators is generated by $k[x,y]$-linear combinations of the columns of
$A_{d-1}$, and since these columns are bigraded, the $k[s,t]$-module
of relations will be bihomogeneous and we can write
$$\Coker A_{d-1} =\bigoplus_{0\leq D,\  1 \leq j} \left(\Coker A_{d-1}\right)_{(D,j)} .$$
Consider the $k[s,t]$-module $\left(\Coker A_{d-1}\right)_{(d,d)}$, the bihomogeneous component of
$\Coker A_{d-1}$ of bidegree $(d,d)$.
It is generated by the images of
$$x y^{d-1} \mathbf{e}_1, x^2 y^{d-2} \mathbf{e}_2, \dots, x^{d-2} y^{2} \mathbf{e}_{d-2}, x^{d-1} y \mathbf{e}_{d-1}$$
and the relations among these generators are given by $k[s,t]$-linear combinations of
$$ y^{d-2} \mathbf{c}_2,  x y^{d-3} \mathbf{c}_3, \dots,  x^{d-3} y \mathbf{c}_{d-1}, x^{d-2} \mathbf{c}_{d}$$
where $\mathbf{c}_1,\dots,\mathbf{c}_{d+1}$ are the columns of
$A_{d-1}$.
So we have $$\left(\Coker A_{d-1}\right)_{(d,d)}= \Coker B_{d-1}$$
where $B_{d-1}$ is viewed as a $k[s,t]$-homomorphism $k[s,t]^{d-1} \rightarrow k[s,t]^{d-1} .$

Using Lemma \ref{lemma1}(i) we deduce that for all $d\geq 2$
the direct summand $\left(\Coker A_{d-1}\right)_{(d,d)}$ of $\Coker A_{d-1}$ has $\tau_{d-1}$ torsion, and
so does $\Coker A_{d-1}$ itself.

Lemma \ref{lemma1}(ii) applied with $\mathcal{S}=\mathbb{N}$
now shows that there exist infinitely many irreducible homogeneous polynomials
$\{ p_i\in k[s,t] : i\geq 1 \}$ each one of them contained in some associated prime of the $R_0$-module
$\oplus_{d\geq 2} \Coker A_{d-1}$.
Clearly, if $i\neq j$ then any prime ideal $P\subset R_0$ which contains
both $p_i$ and $p_j$ must contain both $s$ and $t$.

Since the localisation of $\left(\Coker A_{d-1}\right)_{(d,d)}$ at $s$ does not vanish,
there exist $P_i, P_j\in \Ass_{R_0} \Coker A_{d-1}$ which do not contain $s$ and such that
$p_i\subset P_i$, $p_j\subset P_j$, and the previous paragraph shows that $P_j\neq P_j$.

The second statement now follows from the fact that  $H^2_{R_+}(R)$ is $R_0$-isomorphic to
$\oplus_{d\geq 2} \Coker A_{d-1}$.
\end{proof}

\begin{cor}
Let $T$ be the localisation of $R$ at the irrelevant maximal ideal
$\mathfrak{m}=\langle s,t,x,y,u,v \rangle$. Then $H^2_{(u, v)T} (T)$ has infinitely many associated primes.
\end{cor}
\begin{proof}
Since $\tau_i\in \mathfrak{m}$ for all $i\geq 1$,
$ H^2_{(u,v)T} (T) \cong (H^2_{(u,v)R} (R))_\mathfrak{m}$
has $\tau_i$-torsion  for all $i\geq 1$.
\end{proof}

\section{\bf A connection with associated primes of Frobenius powers}

In this section we apply a technique similar to the one used in section 1 to give a proof of
a slightly more general statement of Theorem 12 in \cite{K}.
The new proof is simpler, open to generalisations and it gives a connection between
associated primes of Frobenius powers of ideals and of local cohomology modules, at least on a purely formal level.

Let $k$ be any field, let $S=k[x,y,s,t]$, let $F=xy(x-y)(sx-ty)=sx^3y-(t+s)x^2y^2+txy^3$ and let
$R=S/FS$.

\begin{thm}\label{theorem2}
Let  $\mathcal{S}$ be an infinite set positive integers and suppose that either $k$ has characteristic zero or that
$k$ has characteristic $p$ and that $\mathcal{S}$ contains infinitely many powers of $p$. The set
$$\bigcup_{n\in \mathcal{S}} \Ass_{R} \left(\frac{R}{\langle x^n, y^n \rangle} \right)$$
is infinite.
\end{thm}
\begin{proof}
We introduce a grading in $S$ by setting $\deg(x)=\deg(y)=1$ and $\deg(s)=\deg(t)=0$. Since $F$ is homogeneous, $R$ is also
graded.

Fix some $n>0$ and consider
the graded $R$-module $T=R/\langle x^n, y^n\rangle$. For each $d>4$
consider $T_d$, the degree $d$ homogeneous component of $T$, as a $k[s,t]$-module.
If $d<n$, $T_d$ is generated by the images of
$y^d, xy^{d-1}, \dots, x^{d-1}y, x^d$
and the relations among these generators are obtained from
$y^{d-4} F, x y^{d-5} F, \dots, x^{d-5}y F, x^{d-4}F$. Using these generators and relations, in the given order,
we write
$T_d=\Coker M_d$ where
$$M_d=\left(
\begin{array}{ccccc}
  0 & 0 & \dots &  & 0 \\
  t &  &  &  &  \\
  -t-s & t &  &  &  \\
  s & -t-s &  &  &  \\
   & s & \ddots & &  \\
   &  &  & t &  \\
   &  &  & -t-s & t \\
   &  &  & s & -t-s \\
   &  &  &  & s \\
  0 & 0 & \dots &  & 0
\end{array}
\right) .$$

When $d=n$, $T_d$ is isomorphic to the cokernel of the submatrix of $M_d$ obtained by deleting the first and last rows
which correspond to the generators $y^n, x^n$ of $T_n$.

When $d=n+1$, $T_d$ is isomorphic to the cokernel of the submatrix of $M_d$ obtained by deleting the first two rows and and last two rows
which correspond to the generators $y^{n+1}, xy^{n}, x^n y, x^{n+1}$ of $T_{n+1}$,
and the resulting submatrix is $B_{n-2}$ defined in Lemma  \ref{lemma1};
the result now follows from that lemma.
\end{proof}

This technique for finding associated primes of non-finitely generated graded modules and of sequences of graded
modules has been applied in \cite{BKS} and \cite{KS} to yield further new and surprising properties of top
local cohomology modules.

\section*{\bf Acknowledgment}

I would like to thank Rodney Sharp, Anurag Singh and Gennady Lyubeznik for reading a first draft of
this paper and for their helpful suggestions.

\end{document}